\documentclass{amsart}

\usepackage{graphicx}
\usepackage{amsmath}
\usepackage[cp1250]{inputenc}
\newcommand\Ker{\mathop{\rm Ker}}
\newcommand\im{\mathop{\rm Im}}

\newtheorem{theorem}{Theorem}[section]

\newtheorem{lemma}[theorem]{Lemma}
\theoremstyle{definition}
\newtheorem{definition}[theorem]{Definition}
\newtheorem{remark}{Remark}
\newtheorem{example}{Example}

\begin{document}

\title[$p$-regularity theory. Applications and developments]{$p$-regularity theory.\\ Applications and developments}

\author{Ewa Bednarczuk$^{2}$, Agnieszka Prusi\'nska $^{1}$, Alexey \ Tret'yakov $^{2,3}$}

\address{1. Department of Mathematics and Physics, Siedlce University, 3-go Maja 54, 08-110 Siedlce, Poland}
\email{aprus@uph.edu.pl}

\address{2. System Research Institute, Polish Academy of Sciences, Newelska 6, 01-447 Warsaw, Poland}

\address{3. Dorodnitsyn Computing Center, Russian Academy of Sciences, Vavilova 40, 119991 Moscow, Russia}



\date{\today}

\begin{abstract}
We present recent advances in the analysis of operator equations with singular operators and constrained optimization problems with constraints given by singular mappings obtained within the framework of the $p$-regularity theory developed over the last twenty years. In particular, we address the problem of description of the tangent cone to the solution set of the operator equations, optimality conditions and solution methods for optimization problems.
\end{abstract}

\maketitle

\section*{Introduction}\label{intro}

We consider problems of solving nonlinear equation
of the form
\begin{equation}\label{eq1}
  F(x)=0,
\end{equation}
and optimization problems of the form
\begin{equation}\label{eq1a}
\min \phi(x) \; \hbox{subject to } \; F(x)=0,
\end{equation}
where $\phi: X\rightarrow \mathbb{R}$ and $F:X\rightarrow Y$ is a sufficiently smooth mapping from a Banach space $X$ to a
Banach space $Y$. Nonlinear problems \eqref{eq1} and \eqref{eq1a} can be divided into two
classes: regular (nonsingular) and singular depending on whether $F$ is regular or singular. Roughly speaking, regular mappings are those
for which implicit function theorem arguments can be applied and singular problems
are those for which they cannot, at least not directly.

In this work, we give an overview of methods and tools of the $p$-regularity
theory in application to the investigation of singular (irregular, degenerate) mappings and singular equality constrained optimization problems. The purpose of this paper is to present selected works in this area in a coherent way, which have been scattered throughout various references.

\section{Essential nonlinearity and singular maps}\label{essential}

Fix a point $x^{\ast}\in X$ and suppose that $F : X \rightarrow Y$ is $\mathcal{C}^1(W)$, where $W$ is a neighborhood of $x^{\ast}$.
The mapping $F$ is \textit{regular} at $x^{\ast}$, if
\begin{equation}\label{eq2}
\im\ F'(x^{\ast})=Y.
\end{equation}
The following lemma on the local representation of regular mapping holds.
\begin{lemma}[Lemma 1., Sec.1.3.3.of \cite{IzTr94}]
	\label{lemma1}
	Let $X$ and $Y$ be Banach spaces, $F : X \rightarrow Y$ and $F\in \mathcal{C}^1(W)$, where $W$ is a neighborhood of $x^{\ast}$.
	
	If $F$ is regular at $x^{*}$, there exist a neighborhood $U$ of $0$ and a neighborhood $V$ of
	$x^{*}$ and a diffeomorphism $\varphi:U\rightarrow V$ such that
	\begin{enumerate}
		\item[(i)] $\varphi(0)=x^{*}$,
		\item[(ii)] $F(\varphi(x))=F(x^{*})+F'(x^{*})x$ for all $x\in U$,
		\item[(iii)]  $\varphi'(0) = I_X$ (the identity map on $X$).
		\end{enumerate}
\end{lemma}
Lemma \ref{lemma1} says that  the diffeomorphism $\varphi$  locally transforms $F$ into the affine mapping,
\begin{equation}\label{eq3}
F(\varphi(x))=F(x^{\ast})+F'(x^{\ast})x \ \ \text{for all } x\in U.
\end{equation}
This fact is also  referred to as the local 'trivialization theorem' (Theorem~1.26 of \cite{ioffe}).
 If the regularity
condition \eqref{eq2} is not satisfied,  the local trivialization of $F$ is not possible ($\varphi$ does not exist), in general.

There exist numerous mappings which do not admit local trivialization.
 The concept of essentially nonlinear mappings defined in \cite{TrMa03}  formalizes this situation.

 \begin{definition}\label{def2.1}
Let $V$ be a neighborhood of $x^{\ast}$ in $X$ and $U\subset X$ be a neighborhood of $0$. A
mapping $F:V\rightarrow Y$, $F\in  \mathcal{C}^2(V)$, is
 \textit{essentially nonlinear at } $x^{\ast}$ if there exists a perturbation of the form
$\widetilde{F}(x^{\ast}+x)=F(x^{\ast}+x)+\omega(x), \hbox{ where } \|\omega(x)\|=o(\|x\|),$
which cannot be trivialized, i.e.  there does not exist any diffeomorphism (i.e. a nondegenerate transformation of coordinates) $\varphi(x):U\rightarrow V$,  $\varphi\in \mathcal{C}^1(U)$,  such that $\varphi(0)=x^{\ast}$, $\varphi'(0)=I_X$ and \eqref{eq3} holds with $\varphi$ and $\widetilde{F}$.
\end{definition}

\begin{definition}
We say the mapping $F$ is \textit{singular} (or  \textit{degenerate})
at $x^{\ast}$ if it fails to be regular; that is, its derivative is not onto:
\begin{equation}\label{eq4}
  \hbox{Im}\ F'(x^{\ast})\neq Y.
\end{equation}
\end{definition}

Let us note that,  if $F$ is singular at the point $x^{\ast}$, $F(x^{*})=0$, i.e.,  there exists $0\neq \xi\in Y$, $\|\xi\|=1$,
\begin{equation}\label{eq4aa}
\xi\not\in \im F'(x^{\ast})
\end{equation}
then $F$ must be essentially nonlinear at $x^{*}$. Indeed,  suppose that $F$ is not essentially nonlinear at $x^{\ast}$ and
define the mapping $\tilde{F}:V\rightarrow Y$ as
\begin{equation}\label{eq4a}
  \tilde{F}(x^{\ast}+x):=F(x^{\ast})+F'(x^{\ast})x+\xi\|x\|^2.
\end{equation}
Note that $\xi\|x\|^2\not\in \im F'(x^{\ast})$ for any $x\in V$.
By Definition~\ref{def2.1}, there exist a neighborhood $U$ of $0$ and a  mapping $\varphi(x):U\rightarrow V$, $\varphi\in \mathcal{C}^1(U)$, such that $\varphi(0)=x^{\ast}$, $\varphi'(0)=I_X$ and
\begin{equation}\label{eq4b}
   \tilde{F}(\varphi(x))=\tilde{F}(x^{\ast})+\tilde{F}'(x^{\ast})x=F(x^{\ast})+F'(x^{\ast})x\ \ \ \text{for  } x\in U.
\end{equation}

 Since  $F'(x^{\ast})x\in \im F'(x^{\ast})$, by \eqref{eq4b} we have
 \begin{equation}\label{eq4c}
  \tilde{F}(\varphi(x))\in \im F'(x^{\ast}).
 \end{equation}
 On the other hand, $\varphi(0)=x^{\ast}$ and $\varphi'(0)=I_X$, and

 \begin{equation}\label{eq4d}
 \begin{split}
    \tilde{F}(\varphi(x)) =& \tilde{F}(x^{\ast}+(\varphi(x)-x^{\ast}))\\
    =& F(x^{\ast})+F'(x^{\ast})(\varphi(x)-x^{\ast})+ \xi\|\varphi(x)-x^{\ast}\|^2\\
    =&  F'(x^{\ast})(\varphi(x)-x^{\ast})+ \xi\|\varphi(0)+ \varphi'(0)x+\omega_1(x)-x^{\ast}\|^2 \\
    =& F'(x^{\ast})(\varphi(x)-x^{\ast})+ \xi\|x+\omega_1(x)\|^2,
 \end{split}
\end{equation}
 where $\|\omega_1(x)\|=o(\|x\|)$. Thus, for small $x$,
 $$\xi\|x+\omega_1(x)\|^2\neq 0.$$
 Taking into account \eqref{eq4aa}, \eqref{eq4d} and the fact that $F'(x^{\ast})(\varphi(x)-x^{\ast})\in \im F'(x^{\ast})$, we conclude from this that
 $$\tilde{F}(\varphi(x))\not\in \im F'(x^{\ast}).$$
 This contradicts \eqref{eq4c} and therefore $F$ is essentially nonlinear at $x^{\ast}$.

The following theorem (see \cite{TrMa03}) establishes the relationship between essential nonlinearity and irregularity.

\begin{theorem}[\cite{TrMa03}]\label{th1}
Suppose $F : V \rightarrow Y$ is $\mathcal{C}^2(V)$, where $V$ is a neighborhood of $x^{\ast}$ in $X$ and $F(x^{\ast})=0$. Then
$F$ is essentially nonlinear at the point $x^{\ast}$ if and only if $F$ is singular at the point
$x^{\ast}$.
\end{theorem}

\section{Examples of singular problems}

\subsection{Description of the solution set. Lyusternik theorem.}\label{LT}

Let $X,Y$ be Banach spaces. Consider the   nonlinear equation \eqref{eq1}

$$  F(x)=0,
$$
where   $F:X\rightarrow Y$,
  $F\in \mathcal{C}^{p+1}(X)$, $p\in \mathbb{N}.$
According to Lyusternik theorem (see \cite{LiSo61}), if $F$
  is regular at $x^{*}$, then  $T_1M(x^{\ast})=\Ker F'(x^{\ast})$, where $T_1M(x^{\ast})$ is the  tangent cone to the set $M(x^{\ast})=\{x\in X: F(x)=F(x^{\ast})=0\}$ at the point $x^{\ast}$. The tangent cone to $M$ at $x^*$ is the collection of all tangent vectors to $M$ at $x^*$ i.e. $h$ is a tangent vector to $M$ at $x^*\in M$ if there exist   $\varepsilon>0$ and a function $r:[0,\varepsilon]\rightarrow X$ with the property that for $t\in[0,\varepsilon]$ we have $x^*+th+r(t)\in M$ and
  $$\lim_{t\rightarrow 0} \frac{\|r(t)\|}{t}=0.$$

  If
  $F$ is singular at the solution point $x^{\ast}$,
  then $T_1M(x^{\ast})\neq \Ker F'(x^{\ast})$.

  For example,  if $F(x)=x_1^2-x_2^2+o(\|x\|^2)$ and $x^{\ast}=0$ then $\Ker F'(0)= \mathbb{R}^2$ and moreover $T_1M(0)=\left\{\left(%
\begin{smallmatrix}
  t \\
  t
\end{smallmatrix}%
\right):t\in \mathbb{R}\right\}\cup \left\{\left(%
\begin{smallmatrix}
  t \\
  -t
\end{smallmatrix}%
\right):t\in \mathbb{R}\right\}$, hence $\Ker F'(0)\neq T_1M(0)$.

The problem of description of the solution sets in more general situations (e.g. general systems of inequalities) is qualitatively approached  by means of metric regularity (\cite{DoLeRo02,Io2000,ioffe})
 and via geometrical derivability (\cite{SeTa06}).

\subsection{Optimality conditions. Lagrange multiplier theorem}\label{sec3.2}

Consider the  optimization problem \eqref{eq1a},
$$
\min \phi(x) \quad \hbox{ subject to } \quad F(x)=0,
$$
where $\phi:X\rightarrow\mathbb{R}$, $\phi\in
\mathcal{C}^{2}(X)$ and  $F:X\rightarrow Y$,
$F\in \mathcal{C}^{p+1}(X)$, $p\in \mathbb{N}.$

Let $x^{*}$ be a solution to \eqref{eq1a}. In the regular case, that is if $\ F'(x^{\ast})\cdot X=Y$, there exists $\lambda^{\ast}\in Y^{\ast}$ such that
$\phi '(x^{\ast})=F'(x^{\ast})^{\ast}\cdot \lambda ^{\ast}.$

Let us consider \eqref{eq1a}, where $X=\mathbb{R}^{3}$, $Y=\mathbb{R}^{2}$, $\phi(x)=x_2^2+x_3$ and
$F(x)=\left(%
                                                                            \begin{array}{c}
                                                                            x_1^2-x_2^2+x_3^2 \\
                                                                            x_1^2-x_2^2+x_3^2+x_2x_3 \\
                                                                            \end{array}%
                                                                            \right)$.
 In this case $x^{\ast}=(0,0,0)^T$ and we can easily obtain that
     $\phi'(x^{\ast})=(0,0,1)^T$,
                                                                    $F'(x^{\ast})=\left(%
\begin{array}{ccc}
  0 & 0 & 0 \\
  0 & 0 & 0 \\
\end{array}%
\right)$. However, it is obvious that
$\varphi'(x^{\ast})\neq F'(x^{\ast})^T\cdot \lambda^*$.

There is a vast literature concerning optimality conditions for general regular (satisfying some constraint qualification condition) constrained optimization problems (see e.g. Chapter~3 of~\cite{BoSh2000}).

\subsection{Newton method for singular equations}\label{sec3.3}

Consider the problem of solving nonlinear equation \eqref{eq1}
where $F:X\rightarrow Y,$
$F\in \mathcal{C}^{p+1}(X),\ p\in \mathbb{N}.$ Let $x^{\ast}$ be a solution
to \eqref{eq1}, i.e. $F(x^{\ast})=0$ and let $F$ be singular at $x^{*}$.

In the finite dimensional case, when $X=\mathbb{R}^n$, $Y=\mathbb{R}^n$ and $F(x)=(f_1(x),\ldots,f_n(x))^T$, singularity of $F$ at $x^{*}$ means that the Jacobian $F'(x^{\ast})$  of $F$ at $x^*$
 is singular as in the following  example.

 \begin{example}[\cite{SzPrTr12}]\label{ex1}

Let $
    F(x)=\left(
           \begin{array}{c}
             x_1+x_2 \\
             x_1x_2 \\
           \end{array}
         \right),
$
$F:\mathbb{R}^2 \rightarrow\mathbb{R}^2,$ where $x^{\ast}=(0,0)^T$ is a
solution to \eqref{eq1} and $
    F'(x^{\ast})=\left(
           \begin{array}{cc}
             1 & 1 \\
             0& 0 \\
           \end{array}
         \right)
$ is singular (degenerate) at $x^{\ast}.$

Let $x_0=(x_{01},x_{02})^T$ and $ x_0
\in U_\varepsilon(0), \ \varepsilon>0 $ be sufficiently
small. Then, for classical Newton method, i.e.
\begin{equation}
\label{7}
    x_{k+1}=x_{k}-\{F'(x_k)\}^{-1}F(x_k), \quad k=0,1,2,3,\ldots ,
\end{equation}
we have
\begin{equation*}
    x_1=\frac{1}{x_{01}-x_{02}}\cdot\left(
                                 -x_{01}x_{02},
                                 x_{01}x_{02}
                             \right)^T.
\end{equation*}
If $x_{01}=x_{02}$ then  $\{F'(x_0)\}^{-1}$ does not exist, hence \eqref{7} is not applicable.

But even ever $ \{F'(x_0)\}^{-1}$ exists, e.g. for
$x_0=(t+t^3,t)^T$, we have $x_1=\left(-\frac{1}{t}-t,
\frac{1}{t}+t \right)^T$ and  $\|x_1-0\|\approx\frac{1}{t}\rightarrow\infty,$ when $t\rightarrow 0.$
For instance, if $t=10^{-5}$ then $\|x_1-0\|\approx10^5$ and we have rejecting
effect.
\end{example}

\begin{example} \cite{Reddien78} Let $F:\mathbb{R}^2 \rightarrow\mathbb{R}^2,$
$$
F(x):=\left(\begin{array}{l}
x_{1}+x_{1}x_{2}+x_{2}^{2}\\
x_{1}^{2}-2x_{1}+x_{2}^{2}
\end{array}\right).
$$

The  singular  root  is $x^{*}=(0,0)^{T}$,
null  space  is
$\Ker F(x^{*})=\text{span}\{(0,1)\}$
and  range
space  is
$\im F(x^{*})=\text{span}\{(1,-2)\}.$
The  Jacobian
$F'(x)$
is  singular  on  the  hyperbole
given by
$$
2x_1-2x_{1}^{2}+6x_{2}-4x_{1}x_{2}+2x_{2}^{2}=0.
$$
\end{example}
For the overview of the existing approaches to Newton-like methods for singular operators, see e.g.
\cite{BKL2010}.

\subsection{Newton method for unconstrained optimization problems}
Consider the following problem,
\begin{equation*}
\min_{x\in\mathbb{R}^{2}}\phi(x)
\end{equation*}
and the scheme
\begin{equation}\label{eq_12}
    x_{k+1}=x_{k}-\{\phi''(x_k)\}^{-1}\phi'(x_k),
\end{equation}
where $ \phi:\mathbb{R}^2 \rightarrow \mathbb{R}$, $ \phi(x)=x_1^2+x_1^2x_2+x_2^4$ (see \cite{SzPrTr12}).

The solution of the considered problem is $x^{\ast}=(0,0)^T$. At the initial point,
$x_0=(x_{01},x_{02})^T$ where $x_{01}=x_{02}\sqrt{6(1+x_{02})}$ we have \\
$$
\phi''(x_0)=\left(%
\begin{array}{cc}
  2+2x_{02} & 2x_{02}\sqrt{6(1+x_{02})} \\
  2x_{02}\sqrt{6(1+x_{02})} & 12x_{02}^2 \\
\end{array}%
\right)
$$ and
$\det \phi''(x_0)=0,$  hence does not exist $\{ \phi''(x_0)\}^{-1}$ and it follows that \eqref{eq_12} is not applicable.

\subsection{Singular problems of calculus of variations}

Consider the following Lagrange problem (see \cite{PrSzTr13}):
\begin{equation}\label{eq13}
J_{0}(x)=\int_{t_1} ^{t_2} f(t,x,x^{\prime}) \mathrm{d}t\rightarrow \min
\end{equation}
subject to the subsidiary conditions
\begin{equation}\label{eq14}
G(x)=\tilde{G}(t,x,x^{\prime})=0, \;\;\;  q(x(t_1),x(t_2))=0
\end{equation}
where $ x\in \mathcal{C}^1([t_1,t_2],\mathbb{R}^n),$
$G:X\rightarrow Y,$ $G\in \mathcal{C}^{p+1}(X),$ $X=\mathcal{C}^1([t_1,t_2],\mathbb{R}^n),$
$Y=\mathcal{C}^1([t_1,t_2],\mathbb{R}^m), $
$\tilde{G}(t,x,x^{\prime})=(G_1(t,x,x^{\prime}),\ldots,G_m(t,x,x^{\prime})),$
$q: \mathbb{R}^n \times \mathbb{R}^n\rightarrow
\mathbb{R}^k .$  We assume that all mappings and their
derivatives are continuous with respect to the corresponding variables
$ t,$ $x,$ $x'$.

In the regular case, if $ \im G^{\prime}(x^{\ast})=Y,$ where
$x^{\ast}(t)$ is a solution to \eqref{eq13}--\eqref{eq14}, then
necessary conditions of Euler-Lagrange
\begin{equation}\label{eq15}
f_x+\lambda(t)G_x-\frac{d}{dt}(f_{x^{\prime}}+\lambda(t)G_{x^{\prime}})=0
\end{equation}
are satisfied.

Let $\lambda:=(\lambda_1,\ldots,\lambda_m)^T,\ $
$\lambda(t)G:=\lambda_1(t)G_1+\cdots+\lambda_m(t)G_m,\ $
$\lambda(t)G_x:=\lambda_1(t)G_{1x}+\cdots+\lambda_m(t)G_{mx}.$
In the singular case, when
$
\im G'(x^{\ast})\neq Y,
$
we can only guarantee that the following equations
\begin{equation}\label{eq17}
\lambda_0 f_x+\lambda(t)G_x-\frac{d}{dt}(\lambda_0f_{x'}+\lambda(t)G_{x'})=0
\end{equation}
are satisfied, where $\lambda_0^2+\|\lambda (t)\|^2=1,$ i.e. $\lambda_0$
might be equal to $0$ and then we have not got any
constructive information on $f.$

\begin{example}[\cite{PrSzTr13}]\label{ex5}
	
	\noindent
	Consider the following problem
	\begin{equation}\label{eq18}
	J_{0}(x)=\int_{0} ^{2\pi} (x_1^2+x_2^2+x_3^2+x_4^2+x_5^2)
	\mathrm{d}t\rightarrow \min
	\end{equation}
	subject to
	\begin{equation}\label{eq19}
	G(x)=\left(%
	\begin{smallmatrix}
	x_1^{\prime}-x_2+x_3^2x_1+x_4^2x_2-x_5^2(x_1+x_2) \\
	x_2^{\prime}+x_1+x_3^2x_2-x_4^2x_1-x_5^2(x_2-x_1) \\
	\end{smallmatrix}%
	\right)=0,
	\end{equation}
	$$x_i(0)=x_i(2\pi), \quad i=1,\ldots,5. $$ Here  $
	f(x):=x_1^2+x_2^2+x_3^2+x_4^2+x_5^2$ and  $q_i(x(0),x(2\pi)):=x_i(0)-x_i(2\pi),$ $i=1,\ldots,5.$
	\newline
	The solution of \eqref{eq18}--\eqref{eq19} is $x^{\ast}(t)=0$ and
	$G^{\prime}(x^{\ast}(t))$ is singular.
	
	Indeed, $
	G^{\prime}(0)=\left(%
	\begin{array}{c}
	(\cdot)^{\prime}_1-(\cdot)_2 \\
	(\cdot)^{\prime}_2+(\cdot)_1 \\
	\end{array}%
	\right)
	$
	and
	$
	G^{\prime}(0)x=\left(%
	\begin{array}{c}
	x^{\prime}_1-x_2 \\
	x^{\prime}_2+x_1 \\
	\end{array}%
	\right).
	$
	
	Let $z(t):=x_1(t)$. Thus, we can consider the following equivalent
	problem: whether the mapping $ \tilde{G}(z)=z^{\prime\prime}+z, \;\;\;z(0)=z(2\pi)$ is surjection or not.
	
	It is obvious that for $y\in \mathcal{C}[0,2\pi],$ such that
	$$ \int_{0} ^{2\pi}\sin \tau\ y(\tau)\mathrm{d}\tau\neq 0 \; \hbox{ or }\; \int_{0} ^{2\pi}\cos \tau\ y(\tau)\mathrm{d}\tau\neq 0,$$ the equation
	$z^{\prime\prime}+z=y$ does not have a solution.
	
	The corresponding Euler-Lagrange equations in
	this case are as follows:
	\begin{eqnarray*}\label{eq20}
		2\lambda_0x_1+\lambda_2-\lambda'_1+\lambda_1x_3^2+\lambda_2 x_5^2-\lambda_2 x_4^2 &=& 0 \\
		2\lambda_0x_2-\lambda_1-\lambda'_2+\lambda_1x_4^2+\lambda_2x_3^2-\lambda_1x_5^2-\lambda_2x_5^2 &=& 0 \\
		2\lambda_0x_3+2\lambda_1x_1x_3+2\lambda_2x_2x_3 &=& 0 \\
		2\lambda_0x_4+2\lambda_1x_2x_4- \lambda_2x_1x_4 &=& 0 \\
		2\lambda_0x_5-2\lambda_1x_5x_1- 2\lambda_1x_2x_5-2\lambda_2x_2x_5+2\lambda_2x_1x_5 &=& 0 \\
		\lambda_i(0)=\lambda_i(2\pi),\;\; i=1,2.& &
	\end{eqnarray*}
	Unfortunately, we cannot guarantee that $\lambda_0\neq 0$ and for
	$\lambda_0 = 0$ we obtain the series of spurious
	solutions to the system \eqref{eq18}--\eqref{eq19}:
	$$ x_1=a\sin t, \; x_2=a\cos t, \;x_3= x_4=x_5=0,\;\lambda_1=b\sin t,\; \lambda_2=b\cos t,\; a,b \in
	\mathbb{R}.$$
\end{example}

\subsection{Modified Lagrange function method}\label{augmented}

Consider the following constrained optimization problem
\begin{equation}\label{9}
\min \phi(x) \quad  \hbox{subject to } \quad g_i(x)\leq 0, \quad i=1,\ldots,m,
\end{equation}
where $\phi:\mathbb{R}^n\rightarrow \mathbb{R}$, $g_i:\mathbb{R}^n\rightarrow \mathbb{R}$
and the modified Lagrangian function $L_E(x,\lambda)$, $L_E:\mathbb{R}^{n+m}\rightarrow \mathbb{R}$ associated with \eqref{9}  (see e.g. \cite{BrEvTr06,Ev77}, cf.  \cite{Be16}),
\begin{equation*}
L_E(x,\lambda):=\phi(x)+\frac{1}{2}\sum_{i=1}^m \lambda_i^2g_i(x).
\end{equation*}
This modification allows to replace a nonlinear optimization problem with a system of nonlinear equations.
Moreover, let us define the mapping $G:\mathbb{R}^n\times \mathbb{R}^m\rightarrow \mathbb{R}^{n+m},$
\begin{equation}\label{gie}
G(x,\lambda):=\left(%
\begin{array}{c}
\nabla \phi(x)+ \frac{1}{2}\sum\limits_{i=1}^m \lambda_i^2 \nabla g_i(x)\\
D(\lambda) g(x) \\
\end{array}%
\right),
\end{equation}
where $D(\lambda):= \hbox{diag}\{\lambda_i\},$ $i=1,\ldots,m$, $\lambda\in \mathbb{R}^m$.

Consider the equation,
\begin{equation}\label{11}
G(x,\lambda)=0_{n+m}.
\end{equation}
For $G(x,\lambda)$  the Jacobian matrix $G'(x,\lambda)$ is given by
\begin{equation*}
G'(x,\lambda)=\left(%
\begin{array}{cc}
\nabla^2 \phi(x)+ \frac{1}{2}\sum\limits_{i=1}^m \lambda_i^2 \nabla^2 g_i(x)&(g'(x))^T D(\lambda) \\
D(\lambda) g(x)& D(g(x)) \\
\end{array}%
\right).
\end{equation*}

If the solution point of \eqref{11} is $(x^{\ast},\lambda^{\ast})$, such that $g_i(x^{\ast})=0$ and
$\lambda_i^{\ast}=0$ then strict complementarity condition (SCQ) defined as $I_0(x^*)=\{j=1,2,\ldots,m:\lambda_j^*=0,g_j(x^*)=0\}\neq \emptyset$ fails. Consequently, $G'(x^{\ast},\lambda^{\ast})$  is a degenerate matrix. The example below illustrates the situation.

\begin{example}\cite{BrEvTr06}
Consider the following problem
\begin{equation}\label{ex_9}
\min_{x\in \mathbb{R}^n} (x_1^2+x_2^2+4 x_1 x_2) \quad  \hbox{subject to } \quad x_1\geq 0,\; x_2\geq 0.
\end{equation}
It is easy to see that $x^*=(0,0)^T$ is the solution to \eqref{ex_9} with the corresponding Lagrange multiplier $\lambda^*=(0,0)^T$.

The modified Lagrange function in this case is
$$L_E=x_1^2+x_2^2+4x_1x_2-\frac{1}{2}\lambda_1^2x_1-\frac{1}{2}\lambda_2^2x_2$$
and the Jacobian matrix $G'(x^{\ast},\lambda^{\ast})$ of
$$G(x,\lambda)=\left(%
\begin{array}{c}
2x_1+4x_2-\frac{1}{2}\lambda_1^2\\
2x_2+4x_1-\frac{1}{2}\lambda_2^2 \\
-\lambda_1x_1\\
-\lambda_2x_2\\
\end{array}%
\right)$$
 is singular.
\end{example}

\section{Elements of $p$-regularity theory}

Let us recall the basic constructions of $p$-regularity theory, whose
basic concepts and main results are described e.g. in \cite{IzTr94, Tr87}.

Suppose that the space $Y$ is decomposed into a direct sum
\begin{equation}\label{wz12}
Y=Y_{1}\oplus \ldots \oplus Y_{p},
\end{equation}
 where
$Y_{1}=\overline{\hbox{Im}\ F'(x^{\ast})},\; Z_{1}=Y.$ Let $Z_{2}$ be
closed complementary subspace to $Y_{1}$ (we assume that such
closed complement exists), and let $P_{Z_{2}}:Y\rightarrow Z_{2}$
be the projection operator onto $Z_{2}$ along $Y_{1}.$ By $Y_{2}$
we mean the closed linear span of the image of the quadratic form
$P_{Z_{2}}F^{(2)}(x^{\ast})[\cdot]^{2}.$ More generally, define
inductively,
$$Y_{i}=\overline{\hbox{span }\hbox{Im}\ P_{Z_{i}} F^{(i)}(x^{\ast})[\cdot]^{i}}\subseteq Z_{i},
\;\;i=2, \ldots , p-1,$$ where $Z_{i}$ is a chosen closed
complementary subspace for $(Y_{1}\oplus \ldots \oplus Y_{i-1})$
with respect to $Y,$ $i=2,\ldots,p$ and $P_{Z_{i}}:Y\rightarrow
Z_{i}$ is the projection operator onto $Z_{i}$ along $(Y_{1}\oplus
\ldots \oplus Y_{i-1})$ with respect to $Y,$ $i=2, \ldots,p.$
Finally, $Y_{p}=Z_{p}.$

The order $p$ is chosen as the minimum number for which
\eqref{wz12} holds. Let us define the following mappings
 $$ F_{i}(x)=P_{Y_{i}}F(x),\;\;\; F_{i}:X\rightarrow Y_{i}
\;\; i=1,\ldots,p,$$ where $\;P_{Y_{i}}:Y\rightarrow Y_{i}\;$ is
the projection operator onto $\;Y_{i}\;$ along $(Y_{1}\oplus
\ldots \oplus Y_{i-1}\oplus Y_{i+1}\oplus \ldots \oplus Y_{p})$
with respect to $Y,$ $i=1,\ldots,p.$

\begin{definition}
The linear operator $\Psi_{p}(h) \in {\mathcal{L}}(X,Y_{1}\oplus \ldots
\oplus Y_{p})$, $h\in X, $ $h\neq 0$
\begin{equation}
\label{p-factor}
\Psi_{p}(h)=F'_{1}(x^{\ast})+
F''_{2}(x^{\ast})h+\cdots + F^{(p)}_{p}(x^{\ast})[h]^{p-1},
\end{equation} is
called the \emph{$p$-factor operator}.
\end{definition}

\begin{example}
	\label{ex_4.2}
	For $p=2$ the formula \eqref{p-factor} takes the form
	\begin{equation}
	\label{p-factor1}
	\Psi_{2}(h)=F'_{1}(x^{\ast})+
	F''_{2}(x^{\ast})h,
	\end{equation}
	where $0\neq h\in X$.
	
Consider the operator $F:\mathbb{R}^{2}\rightarrow \mathbb{R}^{2}$, from the Example~\ref{ex1}, where
    $$F(x)=\left(
           \begin{array}{c}
             x_1+x_2 \\
             x_1x_2 \\
           \end{array}
         \right).
$$
It was shown that the Jacobian of $F(x)$ is singular at $x^{\ast}=(0,0)^T$, hence
$\im F'(x^*)=\text{span} \left\{ (1,0)\right\}\neq \mathbb{R}^2$ and hence
   $Y_1=\text{span} \left\{ (1,0)\right\}$ and
 $Y_2=\text{span} \left\{ (0,1)\right\}.$

 To construct $2$-factor operator we use the projections
 \begin{center}
 $P_{Y_1}=\left(
 \begin{array}{cc}
 	1 & 0 \\
 	0 & 0 \\
 \end{array}
 \right) \quad \ $ $ P_{Y_2}=\left(
 \begin{array}{cc}
 0 & 0 \\
 0 & 1 \\
 \end{array}
 \right)$
  \end{center}
 and define the operators $F_1:\mathbb{R}^2\rightarrow Y_1$ and $F_2:\mathbb{R}^2\rightarrow Y_2$. They are as follows,\\
 \begin{center}
 $F_1(x):=\left(
           \begin{array}{c}
             x_1+x_2 \\
             0 \\
           \end{array}
         \right),\quad \ $     $F_2(x):=\left(
           \begin{array}{c}
             0 \\
             x_1x_2 \\
           \end{array}
         \right), $
 \end{center}
Hence, for $h\in \mathbb{R}^2$, the $2$-factor operator has the form
$$\Psi_2(h)(x):=\left(
           \begin{array}{c}
             x_1+x_2 \\
             h_2x_1+h_1x_2 \\
           \end{array}
         \right),\qquad \text{where}\quad  h=(h_1,h_2)^T.$$
It is easy to see that if $h_1\neq h_2$ then $2$-factor operator is surjective.
\end{example}

\begin{definition}\label{def4.3}
 We say that the mapping $F$ is \emph{$p$-regular at $x^{\ast}$} along an element $h$, if
$$\hbox{Im}\ \Psi_{p}(h)=Y.$$
\end{definition}
As we see from the Example~\ref{ex_4.2} a given mapping $F$ may not be  regular with respect to
all $0\neq h\in X$.
\begin{remark}\label{rem1}
The condition of $p$-regularity of the mapping $F$ at the point $x^{\ast}$ along $h$ is equivalent to the following condition
$$\im F^{(p)}_p(x^{\ast})[h]^{p-1}\circ \Ker \Psi_{p-1}(h)=Y_p,$$
where $\Psi_{p-1}(h):=F_1'(x^{\ast})+F_2''(x^{\ast})+\cdots +F_{p-1}^{(p-1)}(x^{\ast})[h]^{p-2}.$
\end{remark}

\begin{definition}
We say that the mapping $F$ is \emph{$p$-regular at $x^{\ast}$} if
 it is $p$-regular along any $ h\in X$ from the set
$$H_p(x^{\ast}):=\{\bigcap_{k=1}^{p}\hbox{Ker}^{k}F_k^{(k)}(x^{\ast})\} \setminus{\{\mathbf{0}\}}\neq \emptyset,$$
where $k$-kernel of the $k$-order mapping $F_k^{(k)}(x^{\ast})$ is
defined as
\begin{equation*}
    \hbox{Ker}^k F_k^{(k)}(x^{\ast}):=\{ \xi \in
    X:F_k^{(k)}(x^{\ast})[\xi]^k=0\}.
\end{equation*}
\end{definition}

In the Example~\ref{ex_4.2} we have $\Ker^1F'_1(x^{\ast})=\text{span} \left\{ \left(1,-1\right)\right\}$ and $\Ker^2F''_2(x^{\ast})=\text{span} \left\{ \left( 1,0 \right)\right\}\cup \text{span} \left\{ \left(0,1 \right)\right\}$. It means that $H_{2}=\emptyset$. As we see, it may happen that $F$ is $p$-regular along some $h\in X$ but $H_{p}=\emptyset$. Hence, according to Definition~\ref{def4.3} $F$ is 2-regular at $x^*$ along any $h\in X$, $h_1\neq h_2$ and is not 2-regular at $x^*$.

For a linear surjective operator  $\; \Psi_{p}(h) :X\mapsto Y$
between Banach spaces we denote by $\; \{\Psi_{p}(h)\}^{-1}\;$ its
\emph{right inverse}. Therefore $\{\Psi_{p}(h)\}^{-1}:Y
\mapsto 2^{X}$ and we have
\begin{equation*}
\{\Psi_{p}(h)\}^{-1}(y)=\left\{x\in X\colon \Psi_{p}(h) x =
y\right\}.
\end{equation*}
We define the \textit{norm} of $\{\Psi_{p}(h)\}^{-1}$ via the
formula
\begin{equation*}
\|\{\Psi_{p}(h)\}^{-1}\|=\sup_{\|y\|=1}\inf \{\|x\|:
x\in\{\Psi_{p}(h)\}^{-1}(y)\}.
\end{equation*}
We say that $\{\Psi_{p}(h)\}^{-1}$ is \emph{bounded}  if
$\|\{\Psi_{p}(h)\}^{-1}\|<\infty.$

\section{Singular problems via $p$-regularity theory}

\subsection{Generalized Lyusternik theorem}

The following theorem gives a description of the solution set in
the singular case (for the proof see \cite{Tr83}).

\begin{theorem}[\cite{Tr83}, Generalized Lyusternik Theorem]\label{th2}  Let $X$ and $Y$ be
Banach spaces and $U$ be a neighborhood of $x^{\ast}\in X.$ Assume
that \(F:X{\rightarrow} Y, \) $F\in C^{p}(U)$ is $p$-regular at
$x^{\ast}.$ Then
\[T_{1}M(x^{\ast})=H_{p}(x^{\ast}).\]
\end{theorem}

The problem of description of the tangent cone to solution set of the operator equation with the singular mappings has been also considered e.g. in \cite{BeTr17,BrTr07,LeSch98,Tr84}.

We now give another version of the Theorem~\ref{th2}.
To state the result, we shall denote by $\hbox{dist}(x,M)$, the \emph{distance function} from
a point $x \in X$ to a set $M$:
\vskip-6pt
$$\hbox{dist}(x,M)=\inf_{y\in M}\|x-y\|, \quad x\in X.$$

\begin{definition}\label{strong}
A mapping $F\in \mathcal{C}^p$ is called \emph{strongly $p$-regular} at a point $x^*$ if there exists $\alpha >0$ such that
$$\sup_{h\in H_{\alpha}} \left\|\{\Psi_p(h)\}^{-1}\right\|<\infty,$$
where
$$H_{\alpha}=\left\{h\in X: \left\|F_i^{(i)}(x^*)[h]^i\right\|\leq \alpha \text{ for all } i=1,\ldots,p, \; \|h\|=1\right\}.$$
\end{definition}

\begin{theorem}[\cite{PrTre16}]\label{th3}
Let $X$ and $Y$ be Banach spaces, and $U$ be a neighborhood of a point $x^{\ast}\in X$. Assume that $F:X\rightarrow Y$ is a $p$-times continuously Fr\'{e}chet differentiable mapping in $U$ and satisfies the condition of strong $p$-regularity at $x^{\ast}$. Then there
exist a neighborhood $U'\subseteq U$ of $x^{\ast}$, a mapping $\xi\mapsto x(\xi):U'\rightarrow X$, and constants $\delta_1>0$ and $\delta_2>0$ such that
$ F(\xi+x(\xi))=F(x^{\ast}),$
\begin{equation}\label{eq7}
  \|x(\xi)\|\leq \delta_1 \sum_{i=1}^p\frac{\|f_i(\xi)-f_i(x^{\ast})\|}{\|\xi-x^{\ast}\|^{i-1}}
\end{equation}
and
$ \|x(\xi)\|\leq \delta_2 \sum_{i=1}^p\|f_i(\xi)-f_i(x^{\ast})\|^{1/i}$ for all $\xi\in U'$.
\end{theorem}
For the proof, see \cite{IzTr94} and \cite{PrTre16}.

\bigskip

Consider the mapping
$ F(x)=
\left(%
\begin{array}{c}
  x_1^2-x_2^2+x_3^2 \\
   x_1^2-x_2^2+x_3^2+x_2x_3 \\
\end{array}%
\right)
$
from the Section~2.2 and recall that $x^{\ast}=(0,0,0)^T.$ It is easy to see that $F'(x^{\ast})=0$,
$$F''(x^{\ast})= \left(
                  \begin{array}{c}
                    \left(\begin{smallmatrix}
 2 & 0& 0 \\
 0 & -2 & 2\\
 0 & 0 & 2
\end{smallmatrix}\right) \\
 \\
\left(\begin{smallmatrix}
 2 & 0& 0 \\
 0 & -2 & 1\\
 0 & 1 & 2
\end{smallmatrix}\right) \\
                  \end{array}
               \right)$$
 and
 $\Ker^2F''(0)=\text{span}\left\{(1,-1,0)\right\}\cup \text{span}\left\{(1,1,0)\right\}.$
The tangent cone at $x^{*}$ in this case is as follows $T_1M(0)=\text{span}\left\{\left(1,-1,0\right)\right\}\cup \text{span}\left\{\left(1,1,0\right)\right\}$. 

 Let $h=(1,1,0)^T$ (or $h=(1,-1,0)^T$) then $\im F''(0)h=\mathbb{R}^2.$ 
It means that the mapping $F(x)$ is $2$-regular at $x^{\ast}=0$ and in this case
$ \Ker^{2}F''(0)=H_2(0)=T_1M(0).$

\subsection{Optimality conditions for $p$-regular optimization
problems}

We define $p$-factor Lagrange function
$$\mathcal{L}_p(x,\lambda,h):=\phi(x)+\left(\sum_{k=1}^p F_k^{(k-1)}(x)[h]^{k-1},\lambda\right),$$
where $\lambda\in Y^{\ast}$ and
$$
  \bar{\mathcal{L}}_p(x,\lambda,h) := \phi(x)+\left(\sum_{k=1}^p \frac{2}{k(k+1)} F_k^{(k-1)}(x)[h]^{k-1},\lambda\right).
$$

To derive optimality conditions for $p$-regular problems we use Definition~\ref{strong}

\begin{theorem}[\cite{Tr84}, Necessary and sufficient conditions for optimality]\label{th4}
Let $X$ and $Y$ be Banach spaces, $\phi \in \mathcal{C}^2(X),$
$F\in \mathcal{C}^{p+1}(X),$ $F:X\rightarrow Y,$
$\phi:X\rightarrow \mathbb{R}.$ Suppose that $h\in
H_p(x^{\ast})$ and $F$ is $p$-regular along $h$ at the point
$x^{\ast}.$ If $x^{\ast}$ is a local solution to the problem
\eqref{eq1a} then there exist multipliers,
$\lambda^{\ast}(h)\in Y^{\ast}$ such that
\begin{equation}\label{wz13}
    \mathcal{L}p'_{x}(x^{\ast},\lambda^{\ast}(h),h)=0.
   \end{equation}
Moreover, if $F$ is strongly $p$-regular at $x^{\ast},$ there
exist $\alpha>0$ and a multiplier $\lambda^{\ast}(h)$ such that
\eqref{wz13} is fulfilled and
\begin{equation}\label{wz14}
    \bar{\mathcal{L}}p_{xx}(x^{\ast},\lambda^{\ast}(h),h)[h]^2\geq \alpha \|h\|^{2}.
\end{equation}
for every $h\in H_p(x^{\ast}),$ then $x^{\ast}$ is a strict local
minimizer to the problem \eqref{eq1a}.
\end{theorem}

\begin{example} Consider the problem from the Section~\ref{sec3.2}
\begin{equation}\label{eq20a}
x_{2}^{2}+x_{3}\rightarrow \min \hbox{ subject to }
F(x)=\left(
\begin{array}{c}
x_{1}^{2}-x_{2}^{2}+x_{3}^{2} \\
x_{1}^{2}-x_{2}^{2}+x_{3}^{2}+x_{2}x_{3}
\end{array}%
\right) =0.%
\end{equation}%
It is easy to verify that the point $x^{\ast }=0$ is a local minimum to
\eqref{eq20a}.

We have shown in the Section~5.1 that $F$ is singular at $x^{*}$ and for $h=(1,1,0)^T$ the mapping $F(x)$
is $2$-regular at $x^{\ast}=0$ along $h.$ Consider the
$2$-factor-Lagrange function with $\lambda_0=1.$ In this case it has the form
\begin{equation*}
    \mathcal{L}_2(x,\lambda(h),h)=x_2^2+x_3+\alpha(x_1-x_2)+\beta(x_1-x_2+x_3),
   \end{equation*}
   where $\lambda(h)=(\lambda_1(h),\lambda_2(h))$ and
   $\lambda_1(h)=(0,0)^T$, $\lambda_2(h)=(\alpha,\beta)^T.$
   Using the equality $\mathcal{L}'_{2\,x}(x^{\ast},\lambda(h),h)=0$
we obtain $\alpha=1$ and $\beta=-1.$ Putting the coefficients into
we have
$\bar{\mathcal{L}}_2(x^{\ast},\lambda(h),h)=\frac{2}{3}x_2^2.$
Therefore, $\bar{\mathcal{L}}''_{2\,xx}(x^{\ast},\lambda(h),h)[h]^2=\frac{4}{3}>0.$
It means that $x^{\ast}$ is a strict local minimizer to
\eqref{eq20a}.
\end{example}

\subsection{$P$-factor Newton method}

Based on the $p$-factor operator
construction we describe a method for solving nonlinear equations of the form \eqref{eq1}, where
$F:\mathbb{R}^n\rightarrow \mathbb{R}^n$ and the matrix $F'(x^{\ast})$ is singular at the solution point
$x^{\ast}$ (see \cite{BrEvTr06,SzPrTr12}).

Let $Y_1=\hbox{Im} F'(x^{\ast}),\;$
$\bar{P}_1=P_{Y_1^{\perp}}$,  $Y_2=\hbox{Im}
\left(F'(x^{\ast})+\bar{P}_1F''(x^{\ast})h\right),$
$\bar{P}_2=P_{Y_2^{\perp}}$,\\
$Y_{k+1}=\im \left(F'(x^{\ast})+ \sum\limits_{i=1}^k
\bar{P}_iF''(x^{\ast})h+\sum\limits_{%
\begin{smallmatrix}
  i_2>i_1 \\
  i_1,i_2\in \{\overline{1,k} \}
\end{smallmatrix}} \bar{P}_{i_2}\bar{P}_{i_1}F^{(3)}(x^{\ast})[h]^2+\cdots+\right.$

 $\left.+\sum\limits_{%
\begin{smallmatrix}
  i_k>\ldots >i_1 \\
  i_1,\ldots ,i_k \in \{\overline{1,k} \}
\end{smallmatrix}}\bar{P}_{i_k}\ldots \bar{P}_{i_1}F^{(k)}(x^{\ast})[h]^{(k-1)}
\right)$ and
$\bar{P}_{k+1}=P_{Y_{k+1}^{\perp}},$ $k=\overline{2,p-1}.$

 Then the principal scheme of the $p$-factor Newton method is as follows
 \begin{equation*}
    x_{k+1}=x_k-\{F'(x_k)+P_1F''(x_k)h+...+P_{p-1}F^{(p)}(x_k)h^{p-1}\}^{-1}\cdot
\end{equation*}
\begin{equation}\label{2.2}
    \cdot (F(x_k)+P_1F'(x_k)h+...+P_{p-1}F^{(p-1)}(x_k)h^{p-1}),
\end{equation}
where $P_1= \sum\limits_{i=1}^{p-1}\bar{P}_i$, $\qquad P_2=\sum\limits_{%
\begin{smallmatrix}
  i_2>i_1 \\
  i_1,i_2\in \{\overline{1,p-1} \}
\end{smallmatrix}}\bar{P}_{i_2}\bar{P}_{i_1},\;  $
$P_{k+1}=\sum\limits_{%
\begin{smallmatrix}
  i_k>...>i_1 \\
  i_1,...,i_k \in \{\overline{1,p-1} \}
\end{smallmatrix}}\bar{P}_{i_k}...\bar{P}_{i_1},$\\ $k=\overline{2,p-1}$
 and $h$ such that   $\|h\|=1$ is fixed.
$P_i, i=\overline{1,p-1}$ are matrices of orthoprojection at the
solution point $x^{\ast}.$ Note that for
    $$F(x^{\ast})+P_1F'(x^{\ast})h+...+P_{p-1}F^{(p-1)}(x^{\ast})h^{p-1}=0$$
the $p$-factor matrix
\begin{equation}\label{2.4}
F'(x^{\ast})+P_1F''(x^{\ast})h+...+P_{p-1}F^{(p)}(x^{\ast})h^{p-1}
\end{equation}
is not singular (it follows from the $p$-regularity along $h$). It means that
$\bar{P}_p=0, $ $Y_p=\mathbb{R}^n.$

Consider the case $p=2$ for the Example~\ref{ex1},
\begin{equation}\label{2.9}
    x_{k+1}=x_k-\{F'(x_k)+P_1F''(x_k)h\}^{-1}\cdot (F(x_k)+P_1F'(x_k)h)
\end{equation}
where $P_1$ is orthoprojection onto $\hbox{Im}
(F'(x^{\ast}))^{\perp}$ and we choose element $h$ $(\|h\|=1)$,  such
that $2$-factor matrix
\begin{equation}\label{2.6}
    F'(x^{\ast})+P_1F''(x^{\ast})h
\end{equation}
is not singular (in fact, it means that $F$ is $2$-regular along $h$).
Then at the solution point the formula
$
    F(x^{\ast})+P_1F'(x^{\ast})h=0
$
is satisfied, hence we can solve the equation\\
$
    F(x)+P_1F'(x)h=0
$
and by virtue of \eqref{2.6}, $x^{\ast}$ is a locally unique
solution.

\begin{theorem}[\cite{SzPrTr12}]\label{th5}
Let $F\in \mathcal{C}^p(\mathbb{R}^n)$ and there exists $h,\|h\|=1$ such
that $p$-factor matrix \eqref{2.4} is not singular. Then for
any $x_0 \in U_\varepsilon(x^{\ast})$ $(\varepsilon>0$
sufficiently small) and for the scheme \eqref{2.2} the inequality
\begin{equation}\label{2.8}
    \|x_{k+1}-x^{\ast}\|\leq c\|x_{k}-x^{\ast}\|^2, k=0,1,2,\ldots
\end{equation}
holds for some constant $c>0$.
\end{theorem}

\begin{example}[\cite{SzPrTr12}]\label{ex5.8} Let
$
    F(x)=\left(
           \begin{array}{c}
             x_1+x_2 \\
             x_1x_2 \\
           \end{array}
         \right),
$
$x^{\ast}=(0,0)^T$. It was shown in the Example~\ref{ex1} that $F$ is singular at $x^{\ast}=(0,0)^T.$
The scheme of $2$-factor Newton method is as follows \eqref{2.9},
where  $P_1=\left(
           \begin{array}{cc}
             0 & 0 \\
             0 & 1 \\
           \end{array}
         \right),$ $h=(1,-1)^T.$
  Then
$$F'(x_k)+P_1F''(x_k)h=\left(
                         \begin{array}{cc}
                           1 & 1 \\
                           x_k^2-1 & x_k^1+1\\
                         \end{array}
                       \right)$$
 and the formula \eqref{2.9} has the form
\begin{equation*}
    x_{k+1}=x_k-\left(
                         \begin{array}{cc}
                           1 & 1 \\
                           x_k^2-1 & x_k^1+1\\
                         \end{array}
                       \right)^{-1} \left(
                                           \begin{array}{c}
                                             x_k^1+x_k^2\\
                                             x_k^1x_k^2+x_k^2-x_k^1  \\
                                           \end{array}
                                         \right)=
\end{equation*}
\begin{equation*}
    =\left(
                         \begin{array}{cc}
                           1 & 1 \\
                           x_k^2-1 & x_k^1+1\\
                         \end{array}
                       \right)^{-1}\left(
                                           \begin{array}{c}
                                             0\\
                                             x_k^1x_k^2
                                           \end{array}
                                         \right).
\end{equation*}
 It means, that
$
     \|x_{k+1}-0\|\leq c\|x_{k}-0\|^2.
$
\end{example}

\bigskip

\begin{example}[\cite{SzPrTr12}] Consider the following problem
$$
 \min_{x\in \mathbb{R}^2} \, \phi(x),
$$
where $\phi:\mathbb{R}^2\rightarrow \mathbb{R}$ is defined by $\phi(x):=x_1^2+x_1^2x_2+x_2^4$. Moreover, let
$
F(x):=\phi'(x)$ where $\phi'(x)=\left(%
\begin{array}{c}
  2x_1+2x_1 x_2 \\
  x_1^2+4x_2^3 \\
\end{array}%
\right),
$
$x^{\ast}=(0,0)^T$.
It is easy to see that $F$ is 3-regular at $x^{\ast}$ along $h=(1,1)^T$ and
\begin{eqnarray*}
  F'(0)+P_1 F''(0)h+P_2 F^{(3)}(0)[h]^2 &=& \phi''(0)+P_1 \phi^{(3)}(0)h+P_2\varphi^{(4)}(0)[h]^2= \\
  &=& \left(%
\begin{array}{cc}
  2 & -11 \\
  2 & 11 \\
\end{array}%
\right)
\end{eqnarray*}
and this matrix is nonsingular.
Here $\bar{P}_1=\left(%
\begin{array}{cc}
  0 & 0 \\
  0 & 1 \\
\end{array}%
\right),$
$\bar{P}_2=\frac{1}{2}\left(%
\begin{array}{cc}
  1 & -1 \\
  -1 & 1 \\
\end{array}%
\right),$ $P_1=\bar{P}_1+\bar{P}_2=\frac{1}{2}\left(%
\begin{array}{cc}
  1 & -1 \\
  -1 & 3 \\
\end{array}%
\right),$
$P_2=\bar{P}_2\bar{P}_1=\frac{1}{2}\left(%
\begin{array}{cc}
  0 & -1 \\
  0 & 1 \\
\end{array}%
\right).$

Consider the 3-factor scheme
\begin{eqnarray*}
  x_{k+1} &=& x_k-\left(\phi''(0)+P_1\phi^{(3)}(0)[h]+P_2\phi^{(4)}(0)[h]^2\right)^{-1}\cdot \\
   & &  \cdot \left(\phi'(x_k)+P_1\phi''(x_k)[h]+P_2\phi^{(3)}(x_k)[h]^2\right).
\end{eqnarray*}

Let us denote $x_k=(x_1,x_2)^T.$ Then

$\|x_{k+1}-0\|=\left\|x_k-\left(%
\begin{array}{cc}
  2 & -11 \\
  2 & 11 \\
\end{array}%
\right)^{-1}\left(%
\begin{array}{c}
  2x_1-11x_2+2x_1 x_2-6x_2^2 \\
  2x_1+11x_2+x_1^2+18x_2^2+4x_2^3 \\
\end{array}%
\right) \right\|=$\\
$=\frac{1}{44}\left\|\left(%
\begin{array}{c}
  11x_1^2+132x_2^2+22x_1 x_2+44x_2^3 \\
  2x_1^2+48x_2^2-4x_1x_2+8x_2^3 \\
\end{array}%
\right) \right\|\leq 10 \|x_{k}-0\|^2.$
\end{example}

\subsection{Optimality conditions for $p$-regular problems of calculus of variations}

To formulate optimality conditions for singular problems of the form \eqref{eq13}--\eqref{eq14}
 we define $p$-factor Euler-Lagrange function
$$S(x) :=f(x)+\lambda(t)G^{(p-1)}(x)[h]^{p-1},$$
where
$G^{(p-1)}(x)[h]^{p-1}:=g_1(x)+g'_2(x)[h]+\ldots +g_p^{(p-1)}(x)[h]^{p-1} $, $\lambda(t)G^{(p-1)}(x)[h]^{p-1}=\left\langle \lambda(t),\left(g_1(x)+g'_2(x)[h]+\ldots+g_p^{(p-1)}(x)[h]^{p-1}\right)\right\rangle$,
$\lambda(t)=(\lambda_1(t),\ldots,\lambda_m(t))^{T}$ and
$g_i(x),$ $i=1,\ldots,p$ are determined for the map $G(x)$ in similar way as $F_i(x),$ $i=1,\ldots,p$ for the mapping $F(x),$ in the Section 4, i.e. $g_k(x)=P_{Y_k}G(x),$ $k=1,\ldots,p.$

Let
$$g_k^{(k-1)}(x)[h]^{k-1}:= \sum_{i+j=k-1}C_{k-1}^i g_{x^i(x')^j}^{(k-1)}(x)[h]^i[h']^{j},\; k=1,\ldots,p,$$
where
$g^{(k-1)}_{x^i(x')^{j}}(x)=g^{(k-1)}_{\underbrace{x\ldots x}_i \underbrace{x'\ldots x'}_j}(x).$

\begin{definition}\label{def5}
We say that the problem \eqref{eq13}--\eqref{eq14} is $p$-regular at $x^{\ast}$ along $h\in \bigcap\limits_{k=1}^{p}\Ker^k g_k^{(k)}(x^{\ast}),$ $\|h\|\neq 0$ if $$\im \left(g'_1(x^{\ast})+\ldots+g_p^{(p)}(x^{\ast})[h]^{p-1}\right)= \mathcal{C}_m[t_1,t_2].$$
\end{definition}

\begin{theorem}[\cite{PrSzTr13}]\label{th6}
Let $x^{\ast}(t)$ be a solution of \eqref{eq13}--\eqref{eq14} and assume that this problem is $p$-regular at $x^{\ast}$ along
$h\in \bigcap\limits_{k=1}^{p}\Ker^k g_k^{(k)}(x^{\ast}).$ Then there exist a multiplier $\hat{\lambda}(t)=(\hat{\lambda}_1(t),\ldots,\hat{\lambda}_m(t))^{T}$ such that the following $p$-factor Euler-Lagrange equation
\begin{eqnarray}\label{wz23}
 & & S_x(x^{\ast}) -\frac{d}{dt}S_{x'}(x^{\ast})=f_x(x^{\ast})+\left\langle\hat{\lambda},\sum_{k=1}^{p}\sum_{i+j=k-1}C_{k-1}^i g_{x^i(x')^j}^{(k-1)}(x^{\ast})h^i(h')^{j}\right\rangle_{x}-\nonumber\\
   &-& \frac{d}{dt}\left[f_{x'}(x^{\ast})+\left\langle \hat{\lambda}(t),\sum_{k=1}^{p}\sum_{i+j=k-1}C_{k-1}^i g_{x^i(x')^j}^{(k-1)}(x^{\ast})h^i(h')^{j}\right\rangle_{x'}\right] =0,\\
   & & \lambda_i(0)=\lambda_i(2\pi),\;\; i=1,2.\nonumber
\end{eqnarray}
holds.
\end{theorem}

The proof of the above theorem is similar to the one for singular isoperimetric problem in \cite{BeTr88} or \cite{KoTr02}.

Consider Example~\ref{ex5}. The mapping $G$ is $2$-regular (it means that in this case $p=2$) at $\bar{x}=(a\sin t, a \cos t, 0,0,0)^T$ along $h=(a\sin t, a \cos t, 1,1,1)^T$.

Consider the following equation
$$f_x(x)+(G'(x)+P_{Y_2}G''(x)h)^{\ast}\lambda=0 $$
which is equivalent to the system of equations
\begin{equation}\label{wz27}
    \left\{
                    \begin{array}{l}
                      2x_1-\lambda'_1+\lambda_2=0 \\
                      2x_2- \lambda'_2-\lambda_1=0\\
                      2 x_3+2\lambda_1a\sin t+2\lambda_2 a \cos t=0\\
                      2x_4+2 \lambda_1 a \cos t-2\lambda_2a \sin t=0\\
                      2 x_5+2\lambda_1 a(\cos t-\sin t)+2\lambda_2 a(\sin t-\cos t)=0.\\
                      \lambda_i(0)=\lambda_i(2\pi), \; i=1,2.
                    \end{array}
                  \right.
\end{equation}

One can verify that the false solutions of \eqref{eq18}--\eqref{eq19}, that is
$$x_1=a\sin t, \; x_2=a\cos t, \;x_3=x_4=x_5=0$$
do not satisfy the system \eqref{wz27} if $a\neq 0.$ It means that $x_1=a\sin t, \; x_2=a\cos t, \;x_3=x_4=x_5$ do not satisfy 2-factor Euler-Lagrange equation \eqref{wz23} from Theorem~\ref{th6}. The only solution to the Example~\ref{ex5} is $x^*=(0,0,0,0,0)^T.$
Indeed, 2-factor Euler-Lagrange equation in this case for $x^*=(0,0,0,0,0)^T$ has the following form
\begin{equation*}
     \left\{
                    \begin{array}{l}
                      -\lambda'_1+\lambda_2=0 \\
                      - \lambda'_2-\lambda_1=0\\
                      2\lambda_1a\sin t+2\lambda_2 a \cos t=0\\
                      2 \lambda_1 a \cos t-2\lambda_2a \sin t=0\\
                      2\lambda_1 a(\cos t-\sin t)+2\lambda_2 a(\sin t-\cos t)=0.\\
                      \lambda_i(0)=\lambda_i(\pi), \; i=1,2,\\
                    \end{array}
                  \right.
\end{equation*}
where the solution is $\lambda_i^{\ast}(t)=0$, $i=1,2.$

\subsection{Modified Lagrange function method for 2-regular problems}

Consider the constrained optimization problem  \eqref{9},
\begin{equation*}
\min \phi(x) \quad  \hbox{subject to } \quad g_i(x)\leq 0, \quad i=1,\ldots,m
\end{equation*}
and the modified Lagrangian function $L_E(x,\lambda)$ defined in Sec.~\ref{augmented},
\begin{equation*}
L_E(x,\lambda):=\phi(x)+\frac{1}{2}\sum_{i=1}^m \lambda_i^2g_i(x).
\end{equation*}
According to Sec.~\ref{augmented}, the matrix
\begin{equation*}
G'(x,\lambda)=\left(%
\begin{array}{cc}
\nabla^2 \phi(x)+\frac{1}{2}\sum_{i=1}^m \lambda_i^2 \nabla^2 g_i(x)& (g'(x))^T D(\lambda) \\
D(\lambda)g(x)& D(g(x)) \\
\end{array}%
\right)
\end{equation*}
is singular at the solution $(x^*,\lambda^*)$ of \eqref{11} such that $g_i(x^{\ast})=0$ and
$\lambda_i^{\ast}=0$.

 We show that the mapping $G(x,\lambda)$ defined by \eqref{gie} is 2-regular at $(x^*,\lambda^*)$.

   Define the set $I(x^{\ast}):=\{j=1,2,\ldots, m: g_j(x^{\ast})=0\}$ of active constraints, the set $I_0(x^{\ast}):=\{j=1,2,\ldots, m: \lambda^{\ast}=0, g_j(x^{\ast})=0\}$ of weakly active constraints, and the set $I_+(x^{\ast}):=I(x^{\ast})\setminus I_0(x^{\ast})$ of strongly active constraints.

 Since $\lambda^{\ast}=0$ and $g_j(x^{\ast})=0$ for all $j\in I_0(x^{\ast}):=\{1,\ldots, s\}$ the rows $(n+1)$ to $(n+s)$ of $G'(x^{\ast},\lambda^{\ast})$ contain only zeros. Define the vector $h\in \mathbb{R}^{n+m}$ as follows
\begin{equation}\label{h}
h^T:=\left(0_n^T,1_s^T,0_{m-s}^T\right)
\end{equation}
and the mapping $\Phi:\mathbb{R}^n\times \mathbb{R}^m$
\begin{equation}\label{eq-fi}
\Phi(x,\lambda):=G(x,\lambda)+G'(x,\lambda)h,
\end{equation}
with $h$ given by \eqref{h}.

The following fact is well known
\begin{lemma}[\cite{BrEvTr06}]\label{lem1}
	Let an $n\times n$ matrix $V$ and $n\times p$ matrix $Q$ be such that $Q$ has linearly independent columns and $\langle Vx,x\rangle>0 \quad \forall_{x\in \Ker Q^T\setminus\{0\}}$.
	Assume moreover that $D_{N}$ is a full rank diagonal $l\times l$ matrix. Then
	$$\bar{A}:=\left(
	\begin{array}{ccc}
	V & Q & 0 \\
	Q^T & 0 & 0 \\
	0 & 0 & D_N \\
	\end{array}
	\right)
	$$
	is a nonsingular matrix.
\end{lemma}

Let $D(\lambda)$ be the diagonal matrix with $\lambda_j$ as the $j$-th diagonal entry. We say that the constraint qualification condition (CQC) is fulfilled if the gradients of active constraints are linearly independent. The second order sufficient optimality condition holds if there exist $\alpha >0$ such that
\begin{equation}\label{cqc}
  z^T\cdot\nabla_{xx}^2L_E(x^{\ast},\lambda^{\ast})z\geq \alpha \|z\|^2
\end{equation}
for all $z\in \mathbb{R}^n$ satisfying the conditions
$$\langle\nabla g_j(x^{\ast}),z\rangle\leq 0, \; j\in I(x^{\ast}).$$

\begin{lemma}[\cite{BrEvTr06}]\label{lem2}
	Let $\varphi, g_i\in \mathcal{C}^3(\mathbb{R}^n)$ ($i=1,\ldots,m$). Assume that the CQC and the second order sufficient optimality conditions are fulfilled at the solution $(x^{\ast},\lambda^{\ast})$ and $\Phi$ is a mapping given by \eqref{eq-fi}. Then the 2-factor operator $ \Phi'(x,\lambda)=G'(x,\lambda)+G''(x,\lambda)h$ is nonsingular at the point $(x^{\ast},\lambda^{\ast})$.
\end{lemma}
This assertion is obtained if in Lemma~\ref{lem1} we set $V=\nabla_{xx}^2L_E(x^{\ast},\lambda^{\ast})$, $D_N=D(g_N(x^{\ast}))$, where $g_N(x):=\left(g_{p+1}(x),\ldots,g_m(x)\right)^T$ and
$$Q=\left[\nabla g_1(x^{\ast}), \ldots, \nabla g_s(x^{\ast}), \lambda^{\ast}_{s+1}\nabla g_{s+1}(x^{\ast}),\ldots,  \lambda^{\ast}_{p}\nabla g_{p}(x^{\ast})\right].$$
Then $\Phi'(x^{\ast},\lambda^{\ast})=\bar{A}$.

Lemma~\eqref{lem2} implies that 2-factor Newton method
\begin{equation}\label{2fm}
w_{k+1}=w_k-\left[G'(w_k)+G''(w_k)h\right]^{-1}\left(G(w_k)+G'(w_k)h\right), \quad k=1,2,\ldots
\end{equation}
can be applied to solve \eqref{11} and we have the following proposition
\begin{theorem}[\cite{BrEvTr06}]
	Let $x^{\ast}$ be a solution to \eqref{9}. Assume that $\varphi, g_i(x)\in \mathcal{C}^3(\mathbb{R}^n)$, $i=1,\ldots,m$, and the constraint regularity condition CRC and the second order sufficient optimality conditions \eqref{cqc} are fulfilled at the point $x^{\ast}$. Then there exists a sufficiently small neighborhood $U_{\varepsilon}(w^{\ast})$ of $w^{\ast}=(x^{\ast},\lambda^{\ast})$ such that  the estimation
	$$\|w_{k+1}-w^{\ast}\|\leq \beta \|w_k-w^{\ast}\|^2,$$
holds for the method \eqref{2fm}, where $w_0\in U_{\varepsilon}(w^{\ast})$ and $\beta>0$ is an independent constant.
\end{theorem}

\section*{Acknowledgments} The results of the research of the first and third author carried out under the research
theme No. 165/00/S were financed from the science grant granted by the Ministry of Science
and Higher Education. \\
The work of the third author was supported also by the Russian Foundation for Basic Research (project No. 17-07-00510, 17-07-00493) and the RAS Presidium Programm (program 27).

\end{document}